\documentclass[12pt]{article}
\usepackage{setspace}
\usepackage{textcomp}
\usepackage{amsmath,amsthm,amssymb}
\usepackage[left=1in,
right=1in,nohead]{geometry}
\usepackage{hyperref}
\usepackage{latexsym}
\usepackage{mathrsfs}
\usepackage{rawfonts}
\usepackage[dvips]{graphicx}
\usepackage{stmaryrd}

\usepackage[sc,osf]{mathpazo}

\def\bct{\begin{center}}
\def\ect{\end{center}}
\def\beg{\begin}
\def\bit{\begin{itemize}}
\def\eit{\end{itemize}}

\def\<{\langle}
\def\>{\rangle}

\def\mbb{\mathbb}
\def\mbbz{\mathbb Z}

\def\mco{\mathcal O}

\def\ni{\noindent}

\def\tn{\textnormal}
\def\wt{\widetilde}

\newtheorem{thm}{Theorem}[section]
\newtheorem{lem}[thm]{Lemma}

\newtheorem{cor}[thm]{Corollary}

\title{Compact Complex Surfaces of Locally Conformally Flat Type} \author{Mustafa Kalafat \and Caner Koca}

\begin{document}
\maketitle
\begin{abstract} We show that if a compact complex surface admits
 a locally conformally flat metric, then it cannot contain a smooth rational curve of odd self-intersection. In particular, the surface has to be minimal. Moreover, we give a list of possibilities of such surfaces.
\end{abstract}


\section{Introduction}

A Riemannian $n$-manifold $(M,g)$ is called {\em locally conformally flat (LCF)} if $M$ has an open cover such that for any open set $U$ of the cover we have a strictly positive smooth function $f:U\to \mbb R^+$ and a diffeomorphism $h: U \longrightarrow \mbb R^n$ such that the pull-back of the Euclidean metric $g_{Euc}$ on $\mathbb R^n$ is conformally related to the restriction of $g$ on $U$; i.e.
$$h^*g_{Euc}=fg.$$

\noindent In this paper, we are specifically interested in dimension four and in the compact case. In particular,
we would like to see which compact complex surfaces can possibly admit an LCF metric.
For this purpose we start with proving the following 
result.

\noindent {\bf Theorem ~\ref{4lcfimpliesselfintersectionsphere0}.} {\em
If a compact complex surface admits a locally conformally flat Riemannian metric, then it cannot contain a smooth rational curve of
odd self-intersection.}

\noindent Since a non-minimal 
complex surface by definition
contains a smooth rational curve $\mbb{CP}_1$ of self-intersection $-1$, we have the following consequence:

\begin{cor} If a compact complex surface admits a locally conformally flat Riemannian metric, then 
it has to be minimal.
\end{cor}

\noindent We apply this corollary to the Enriques--Kodaira classification of surfaces (\cite{bpv}, p.244), and eliminate some of the surfaces appearing
 on the list (see Lemma \ref{listtheorem1}). 
We also analyze the case of elliptic fibrations separately in Theorem \ref{ellipticcase}.
As a consequence of these results, we obtain the following 
 list of possibilities:


\begin{thm}\label{listtheorem2} If a compact complex surface admits a locally conformally flat Riemannian metric, then it must be one of the following surfaces:
\beg{enumerate}
\item\label{ourVII} a Hopf surface, or an Inoue surface with vanishing second Betti number,
\item\label{ourruled} a minimal ruled surface fibered over a curve
$\Sigma_g$
of genus $g\geq 2$,
\item\label{ourelliptic} a minimal elliptic fibration with no singular, 
but possibly with multiple fibers over a genus $g\geq 1$ curve, 
\item a minimal torus which is not elliptic,
\item a non-simply-connected minimal surface of general type of Euler characteristic $\chi\geq 4$ which does not admit a Bergmann metric.
\end{enumerate}
\end{thm}
\ni The spaces that admit a Bergmann metric are of the form $\mathbb C\mathcal H_2/G$, i.e. holomorphic quotient of the complex hyperbolic plane. The interested reader may wish to consult the references
\cite{borel,mumford,chen} for examples of these quotients.

For the final case we conjecture that indeed no surface of general type
admits an LCF metric. One of the intuitions behind this conjecture is that these surfaces have a large
fundamental group, and thus, it seems unlikely that they can be mapped into the group conformal
transformations through the holonomy representation, see section \S
\ref{secdevelopingmap} for the background. The first author obtained a partial result in
this direction in \cite{lcfgeneraltype}. Namely, for product surfaces of general
type, if the holonomy representation is discrete and faithful, then there exists
no LCF metric. On the other hand, holonomy representation can be non-discrete for
this type of metrics. This is contrary to the hyperbolic manifold case. It is a
difficult task to handle the non-discrete representations.

The {\em Weyl invariant} of a compact smooth $n$-manifold $M$ is defined
by $$W(M):=\inf_{g \in {Met}(M)}{\int_M |W_g|^{n/2} \,d\mu_g}$$
where $W_g$ is the Weyl curvature tensor of the metric $g$.
If $n\geq 3$, then any LCF metric $g$ has $W_g\equiv 0$ \cite{besse};  therefore, this invariant turns out to be zero for manifolds with an LCF metric.
If $M$ is a compact quotient of the complex hyperbolic space, then its natural
Bergmann metric attains the minimum by \cite{botvinnikosamuharish} using
a result of LeBrun \cite{clmostow}. This implies that the signature is strictly positive.
Consequently, they obtain $W(M)=48\pi^2\tau>0$, and this prevents the possibility for these surfaces to admit an LCF metric.

On the other hand,
the \emph{Weyl energy} of a product metric $g_{\tn{Prod}}$ on the product of curves $\Sigma_g\times\Sigma_h$ of genera $g$ and $h$ can be computed as
$$W(g_{\tn{Prod}}) :=\int\limits_{\Sigma_g\times\Sigma_h} |W_{g_{\tn{Prod}}}|^2 d\mu =\frac{128\pi^2}{3}(1-g)(1-h)+\frac{2}{3}\int\limits_{\Sigma_g\times\Sigma_h}(\kappa_g-\kappa_h)^2\,d\mu, $$ where $\kappa_g,\kappa_h$ are the Gauss curvatures of each factor \cite{osamukobayashi}. This implies that, for $g,h\geq 2$, the standard product metric (for which $\kappa_g=\kappa_h=$ const.), which is K\"ahler-Einstein, has the minimum Weyl energy among \emph{all product metrics}. Note that the Weyl energy of this K\"ahler-Einstein metric is strictly positive.
However, currently it is not known whether the Weyl energy goes below this level for other (non-product) metrics on $\Sigma_g\times\Sigma_h$, $g,h\geq 2$.

A final remark about the Weyl invariant is that $W(\Sigma_1 \times \Sigma_g)=0$ for any genus $g$. This was first observed by Kobayashi in \cite{osamukobayashi}, as a consequence of his result
which states that the Weyl invariant is zero for manifolds with a free and differentiable circle action. We know that $\Sigma_1\times \Sigma_g$ admits an LCF metric (the flat metric) when $g=1$, but it does not admit an LCF metric when $g=0$.
For higher genera $g\geq 2$, even though the Weyl invariant $W(\Sigma_1\times\Sigma_g)$ is zero (because of the existence of an $S^1$-action induced from the one on the first factor), it is not known whether or not this manifold admits an LCF metric.




The outline of the paper is as follows: In Section \S\ref{secdevelopingmap} we recall the developing map construction for LCF manifolds, and prove the first main result. In \S\ref{seckodairaenriquesclassification} we
obtain a list by analyzing the Kodaira-Enriques classification of
compact complex surfaces. In Section \S\ref{secelliptic} we deal with the elliptic fibration case separately. In \S\ref{secconverse} we give information about the converse case. Finally in \S\ref{sechermitianLCF} we relate our classification to that for the Hermitian case.


\noindent{\bf Acknowledgements.} The authors would like to thank Claude LeBrun for his suggestions and encouragement. The authors also thank Y. G\"urta\c s, A. Akhmedov and Weiwei Wu for useful discussions, and the referee for useful remarks. This work is partially supported by the grant $\sharp$113F159 of T\"UB\.ITAK\footnote{Turkish Science and Research Council.}.


\section{Developing map of a locally conformally flat manifold}\label{secdevelopingmap}

Before defining the developing map of a locally conformally flat manifold,
let us offer some motivation for the definition. We would like to modify the definition of local conformal flatness so that one uses charts and transition maps rather than the Riemannian metric directly. The key theorem in this case is due to Liouville and Gehring (see \cite{louisvillegehring} p.389 or for a recent survey  \cite{ralphhoward}), which states that for $n\geq 3 $ and any open set $U\subset\mbb R^n$, any $C^1$ conformal map $\varphi:U\to \mbb R^n$ is the restriction of a M\"obius transformation of $S^n$.  {\em M\"obius transformations} $\tn{M\"ob}(S^n)=\tn{Conf}(S^n)$ is the group of conformal diffeomorphisms of the round $n$-sphere, and they are generated by inversions in round spheres. So, they constitute a group of real analytic diffeomorphisms of the real analytic manifold $S^n$ by the Liouville-Gehring theorem.
Alternatively, they are the restrictions of the full group of isometries of the hyperbolic space $\mbb L^{n+1}$ to its ideal boundary $S^n$, described as follows. Consider $\mbb R^{n+2}$ with its Lorentzian metric $g_1=dx_1^2+\cdots+dx_{n+1}^2-dx_{n+2}^2$. Let $O(n+1,1)$ be the group of linear maps that preserve the Lorentzian metric. We embed the two mentioned spaces into $\mbb R^{n+2}$ as follows:
$$\mbb L^{n+1}=\{x \in \mbb R^{n+2} : |x|_1^2=-1~~\tn{and}~~x_{n+2}>0\}$$
$$S^n=\{x \in \mbb R^{n+2} : |x|_1^2=0~~\tn{and}~~x_{n+2}=1\},$$
i.e. $\mbb L$ is the upper part of the hyperboloid asymptotic to the light cone and $S^n$ is the unit sphere in the upper light cone which is the boundary of the Klein model $\mbb K$ of the hyperbolic space, see the Figure 5 in \cite{flavorsofgeometrycannon}. The restriction of the Lorentzian metric on $\mbb L^{n+1}$ and $S^n$ gives hyperbolic and round metrics which are positive definite and of constant curvature $-1$ and $1$, respectively. Consider
$$\tn{Isom}(\mbb L^{n+1})=O^+(n+1,1):=\{ A \in O(n+1,1) : A~~\tn{preserves}~~\mbb L \}.$$
We define an isomorphism,
$$\Psi : \tn{Isom}(\mbb L^{n+1}) ~\tilde{\longrightarrow}~ \tn{M\"ob}(S^n),
~~a\mapsto \Psi_a$$
by the following procedure. Take $a\in O^+(n+1,1)$ so that
for $y=(y_1\cdots y_{n+1})$ and
$a(y,1)=(a_1y, a_2y)\in \mbb R^{n+1}\times \mbb R^+$ 
i.e. $a_2:=\pi_{n+2}\circ a \circ \pi_{1\cdots n+1}$
define
$$\Psi_a:S^n\to S^n ~~~\tn{by}~~~ \Psi_a(y,1):=\left( {a_1y \over a_2y},1 \right)$$
This is a conformal map on the sphere since it is the map $y\mapsto (a_1y,a_2y)$, which is an isometry of the sphere on its image, followed by rescaling via the factor  $(a_2y)^{-1}$. So whenever we define a locally conformally flat structure, instead of local conformal diffeomorphisms into $\mbb R^n$, we map into $S^n$. 
\beg{defn} A {\em locally conformally flat} structure on a smooth manifold $M$ is a smooth atlas $\{(U_i,h_i)_{i\in I}\}$ where the maps
$h_i:U_i\to S^n$ are diffeomorphisms onto their images and
the transition maps $h_i\circ h_j^{-1}\in\tn{M\"ob}(S^n)$ after restriction. \end{defn}
Now 
start with one of the flattening (or rounding) maps $h_1:U_1\to S^n$.
Let $\alpha$ be a path in $M$ beginning in $U_1$. We would like to analytically continue $h_1$ along this path. Proceeding inductively,
on a component of $\alpha\cap U_i$ the analytic continuation of $h_1$
is a shift away from $h_i$, i.e. of the form $\Gamma\circ h_i$ for some $\Gamma\in \tn{M\"ob}(S^n)$. This way $h_1$ is analytically continued
along every path of $M$ starting at a point in $U_1$. Therefore, there is
a global analytic continuation $D$ of $h_1$ defined on the universal cover
$\wt{M}$ 
since it is defined as a quotient space of paths in $M$. $D$ is called the
{\em developing map} of the locally conformally flat space.
$$\begin{array}{r@{}ccc}
   &\wt{M}     & \stackrel{D}{\longrightarrow} & S^n \\
 p &\downarrow &                 &     \\
   &M          &                 &
\end{array}$$
If one starts with a different flattening open subset instead of $U_1$,
one gets another developing map which differs from $D$ by a composition with
 a M\" obius transformation. Hence, the developing map is defined uniquely up to a composition with an element in $\tn{M\"ob}(S^n)$.
 This uniqueness property has the following consequence. Given any
 covering transformation $T$ of the universal covering, there is a
unique element $g\in \tn{M\"ob}(S^n)$ such that
$$D\circ T= g\circ D.$$
This correspondence defines a homomorphism
$$\rho: \pi_1(M) \longrightarrow \tn{M\"ob}(S^n)$$
called the {\em holonomy representation} of $M$.
Conversely, starting with a pair $(D,\rho)$ where
$\rho$ is a representation of the fundamental group into M\" obius
transformations and $D$ is any $\rho$-equivariant local diffeomorphism
of $\wt{M}$ into $S^n$, one can construct the corresponding LCF structure
on $M$ by pulling back the standart LCF structure from $S^n$ to $\wt{M}$
via $D$, and then projecting it down.

\begin{thm}\label{4lcfimpliesselfintersectionsphere0}
If a compact complex surface admits a locally conformally flat Riemannian metric, then it cannot contain a rational curve of
odd self-intersection.
\end{thm}
\beg{proof} 
Let $f:S^2\to M$ be a smoothly embedded complex sphere in a compact complex surface $M$. 
Since the fundamental group of the sphere is trivial we have $f_*\pi_1(S^2)\subset p_*\pi_1(\wt{M})$. So by the general lifting lemma (\cite{munkres} p.478), we can lift the embedding to a continuous map $\tilde{f}:S^2\to\wt{M}$ into the universal cover, at any chosen base point in a unique way. Since $p$ is a local diffeomorphism and $f$ is an embedding, the map $\tilde{f}$ is also an embedding locally, hence an immersion. We can conclude that the self-intersection numbers in $M$ and
the universal cover
$$I(f,f)=I(\tilde{f},\tilde{f})$$
are the same since there is a local diffeomorphism and the intersection numbers can be computed through the local deformations of the submanifolds. To be precise, self-intersection number is obtained
by perturbing a copy of the sphere in a neighborhood to make it transverse to itself and counting the signed number of points according to the orientation. Since the covering map is a local diffeomorphism,
it becomes a bijection when restricted to a lifting of the sphere.
Since at the same time it is a local diffeomorphism in a neighborhood of a point, by compactness, passing to a finite cover one can introduce a metric and find a uniform
$\epsilon$ neighborhood on which the covering map is a diffeomorphism.
If the perturbed sphere goes beyond this neighborhood, then we just push it inside without changing the intersection points.

As the second step, we note that the lifted sphere is also a holomorphic one. So that the adjunction formula \cite{donaldsonandkronheimer}
\begin{equation}\label{adjunction}
2g(C)-2 = [C]^2 - c_1(S)[C]
\end{equation}
for a smooth connected curve $C$ of a complex surface $S$ is applicable.
Since the developing map is obtained through local
flattening conformal diffeomorphisms, it is an immersion.
The formula (\ref{adjunction}) has no analogue in the image because
the Chern class is not defined. On the other hand the Stiefel-Whitney class \emph{is} defined. Since $w_2(S^4)=0$, by naturality of characteristic classes we have
$$w_2(T\wt M)=w_2(D^*TS^4)=D^*w_2(TS^4)=0.$$
If one takes the (mod 2) reduction of both sides of (\ref{adjunction}) applied to $C=\tilde f(S^2)$ and $S=\widetilde{M}$, and inserting $c_1(\widetilde M)\equiv w_2 (\widetilde M)\equiv 0 \textnormal{\ (mod 2)}$, one gets
$$0  =  [\tilde{f}(S^2)]^2 ~~\tn{(mod 2)}.$$
\end{proof}

\ni In particular, there cannot be a ($-1$)-self-intersecting smooth rational curve in $M$, and thus, $M$ must be minimal.

\beg{rmk} There are actually immersions $T^*S^2\to \mbb R^4$ which
realize the sphere as a Lagrangian submanifold with respect to the standart symplectic structures of both sides. Note the parity of the self-intersection of the sphere. See \cite{aurouxSLfibrations} 
for details.
\end{rmk}


\section{Kodaira--Enriques classification} 
\label{seckodairaenriquesclassification}

In this section we give a list  of complex surfaces which may possibly admit a locally conformally flat metric.
The idea is to go through the classes of surfaces in the Kodaira--Enriques classification. 
According to the 
classification (\cite{bpv}, p.244), the following is the complete list of minimal surfaces:

\begin{enumerate}
\item Minimal rational surfaces
\item Minimal surfaces of class VII
\item Ruled surfaces of genus $g\geq 1$
\item Enriques surfaces
\item Bi-elliptic surfaces
\item Primary or secondary Kodaira surfaces
\item K3-surfaces
\item Tori
\item Minimal properly elliptic surfaces
\item Minimal surfaces of general type
\end{enumerate}

The assumption that the surface admit an LCF metric helps us to eliminate some of these possibilities by close inspection.
First of all, we make the following general remark: A compact complex surface admitting an LCF metric has to be \textit{of signature $\tau=0$} and \textit{non-simply-connected}.
The fact that the signature $\tau$ is zero follows from the Hirzebruch Signature formula.


\begin{thm}[Hirzebruch] Let  $(M,g)$ be an oriented Riemannian manifold. Then the signature of the manşfold can be expressed in terms of curvature quantities as follows. 
$$\tau = {1\over 12\pi^2}\int_M |W^+|^2 - |W^-|^2 \,d\mu.$$
\end{thm}

\ni Here $W^\pm$ are the self-dual and anti-self-dual parts of the Weyl tensor.
This formula is a combination of two results. One of them is the signature theorem
of Hirzebruch which expresses the signature of an oriented 4-manifold 
as a multiple of the integral of its first Pontrjagin class over the manifold.  See \cite{hirzebruchsignaturetheorem} for a reference. Another result is that using Chern-Weil theory one can express the characteristic classes using curvature quantities. See \cite{chernwithoutpotential}.
Since $W=0$ for any LCF metric \cite{besse}, we see that $\tau$ has to be zero.

Furthermore, we will make use of the following theorems of Kuiper.

\beg{thm}[\cite{kuiper}] \label{kuiperthm} Let $(M^n,g)$ be a simply connected, LCF $n$-manifold of class $C^1$. Then there is a conformal immersion $f:M\to S^n$. If in addition $M$ is compact, then this map is a conformal diffeomorphism.
\end{thm}

\beg{thm}[\cite{kuiper2}] \label{kuiperthm2} Universal cover of a compact, LCF space with an infinite Abelian fundamental group must be $\mbb R^n$ or $\mbb R \times S^{n-1}$.
\end{thm}

\ni According to the first theorem, the 4-sphere $S^4$ is the only compact, simply-
connected 4-manifold with an LCF metric. Since $S^4$ is not a complex manifold, a
compact complex surface with an LCF metric cannot be simply connected.
Now, let us analyze the above list.

The first case is a minimal rational surface. A surface is called \emph{rational} if and only if it is birationally equivalent to
the complex projective plane. 
 The possibilities for the minimal models are
the complex projective plane $\mbb{CP}_2$ and the
Hirzebruch surfaces $\mbb F_n=\mbb P(\mco \oplus \mco_{\mbb{CP}_1} (n))$ for $n=0,2,3 \dots$ (see \cite{bpv}).
The Hirzebruch surfaces fall into two distinct smooth topological types
$S^2 \times S^2$ and $\mbb{CP}_2\sharp\,\overline{\mbb{CP}}_2$ determined by parity of $n$ (see \cite{hirz}). Both of these types are simply-connected, as is $\mbb{CP}_2$. Thus, they cannot admit an LCF metric by Kuiper's theorem.

The second item in the list is the minimal surfaces of class VII. These surfaces are characterized by their Kodaira dimension $\kappa=-\infty$ and Betti number $b_1=1$ (therefore, they are not simply-connected).
Furthermore their Chern numbers satisfy $c_1^2\leq 0$ and $c_2\geq 0$.
Combining with the identity $$c_1^2=2\chi+3\tau=2\chi$$
for LCF complex surfaces, we reach the conclusion that $\chi=0$.
Since $b_1=1$, we can compute the second Betti number as follows, $$0=\chi=2-2b_1+b_2=b_2.$$
Class VII minimal surfaces of vanishing second Betti number are classified by Bogomolov in \cite{bogomolov1,bogomolov2}: Hopf surfaces and Inoue surfaces are the only two possibilities. A surface is called a {\em Hopf surface} if its universal cover is biholomorphic to $\mbb C^2-0$. The other possibility are {\em Inoue surfaces} with $b_2=0$. Their universal cover is biholomorphic to $\mbb C \times \mbb H$, i.e. complex line times the hyperbolic disk.

The third is the case of ruled surfaces of genus $g\geq 1$. Such a surface
admits a ruling, i.e.  a locally trivial, holomorphic fibration over a smooth non-rational curve with fiber $\mbb{CP}_1$ and structural group $\tn{PGL}(2,\mbb C)$. This can be thought of as a projectivization of a complex rank 2-bundle over a Riemann surface.
Now, we claim that the base cannot be a torus: suppose that the base \emph{is} a torus.
Then, topologically, we have the following fiber bundle
$$S^2\longrightarrow M \longrightarrow T^2.$$
The homotopy exact sequence for this bundle involves the following terms:
$$\cdots\to\pi_3 T^2 \to \pi_2 S^2 \to \pi_2 M \to \pi_2 T^2 \to\cdots$$
Here, the terms at the two ends are zero since the universal cover of torus is
contractible. Therefore, we have the isomorphism $\pi_2 M \approx \pi_2 S^2
\approx \mbbz$. Thus, the second homotopy group of the universal cover
$\wt{M}$ is non-trivial. Taking a look at the remaining terms of the homotopy exact sequence on the right we have
$$\cdots\to\pi_1 S^2 \to \pi_1 M \to \pi_1 T^2 \to \pi_0 S^2\to\cdots$$
Again the end terms vanish, and we have the isomorphism
$\pi_1 M \approx \pi_1 T^2 \approx \mbbz\oplus\mbbz$.
This is an infinite abelian group.
However, due to 
the second theorem of Kuiper that we stated above, since the fundamental group is
infinite abelian, the universal cover 
must be $\mbb R^4$ or $\mbb R \times S^3$ if there is any LCF metric.
Ours have non-trivial second homotopy group, so it is none of these. Hence the genus $g=1$ case yields a contradiction. 

The fourth and seventh possibilities are eliminated, because the signatures of Enriques and K3 surfaces are nonzero: $\tau(E)=-8$ and $\tau(K3)=-16$.

Finally, let us consider surfaces of general type. We know that the Chern numbers $c_1^2$ and $c_2$ are strictly positive for these surfaces.
Recall the formula for the holomorphic Euler characteristic:
$12\chi_h=c_1^2+c_2$. Since this is a non-zero positive integer multiple of $12$, we have $c_1^2+c_2\geq 12$.
Recall the identity $c_1^2=2c_2+3\tau$
for complex surfaces. Adding $c_2$ to both sides
and applying the previous inequality we obtain $c_2+\tau\geq 4$.
 Since the signature $\tau$ is zero for LCF surfaces, we get $c_2=\chi\geq 4$.

Now, we can list the remaining cases as follows.

\begin{lem}\label{listtheorem1} If a compact complex surface $(M,J)$ admits a locally conformally flat Riemannian metric, then it can be either
\beg{enumerate}
\item\label{itemhopfinoue} a Hopf surface, or an Inoue surface with $b_2=0$,
\item  a ruled surface fibered over a Riemann surface
$\Sigma_g$
of genus $g\geq 2$,
\item\label{elliptic1} a bi-elliptic surface,
\item a primary or secondary Kodaira surface,
\item a torus,
\item\label{elliptic4} a minimal properly elliptic surface, or
\item a non-simply-connected minimal surface of general type of Euler characteristic $\chi\geq 4$.
\end{enumerate}
\end{lem}


\section{Elliptic Surfaces}\label{secelliptic}

In the case of elliptic surfaces one can actually make a more refined classification. We start with the following classification theorem stated in \cite{gs} p.314, a summary of research done by various people. See the references therein.

\beg{thm}\label{ellipticclassification}
A relatively minimal elliptic surface with nonzero Euler characteristic is diffeomorphic to $E(n,g)_{p_1\cdots p_k}$ for exactly one choice of the integers involved for $$1\leq n,~ 0\leq g,k,~
2\leq p_1\cdots\leq p_k~~\tn{and}~~k\neq 1~~\tn{if}~~(n,g)=(1,0).$$
\end{thm}
\ni Here, {\em relatively minimal} means that the fibers do not contain any sphere of self intersection $-1$. This is a generalization of being minimal. $E(1)$ is defined to be the surface $\mbb{CP}_2\sharp\,9\,\overline{\mbb{CP}}_2$ considered with its elliptic fibration. Then taking its fiber sum with itself $n$-times, one gets
$E(n)$. Furthermore taking the fiber sum with the trivial fibration
$\Sigma_1\times\Sigma_g$  over the Riemann surface of genus $g$, one gets the space $E(n,g)$. Finally the subindices $p_i$ denotes the multiplicity
of a logarithmic transformation. Log transform is a standard way to introduce a
multiple fiber.
Using this classification theorem we can prove our result.

\beg{thm}\label{ellipticcase} If an elliptic surface admits a locally conformally flat
metric then it is minimal, and it has vanishing Euler characteristic and signature. Moreover it has to be a torus bundle over a curve, 
outside the multiple fibers.
\end{thm}
\beg{proof} The signature of $E(n,g)_{p_1\cdots p_k}$ is computed to be
$\tau=-8n$, see \cite{gs}. If we assume that there exists an LCF metric,
then signature has to vanish and $n=0$. Applying the Theorem \ref{ellipticclassification} we reach  the conclusion that the Euler characteristic $\chi=0$. But this is the Euler characteristics of the fiber bundle: $\chi=\chi(\tn{fiber})\times \chi(\tn{base}) =0$.
Since a cusp fiber contributes by $2$ and a fishtail fiber contributes by $1$ to the Euler characteristic, there are no singular fibers.
\end{proof}
\ni {\em Logarithmic tranformation} is a standard way to introduce a multiple fiber. Topologically, picking up a lattitute $l$ of a smooth fiber, multiplying with the disc in the base, replacing the solid torus by another solid torus which has multiple Seifert fibered central circle is basically what this operation means. Note that this does not change the Euler characteristic. It changes the homology class of the fiber.

An elliptic surface fibered over a rational curve is either a product or a Hopf surface, see \cite{bpv} p.196. Since the product does not admit LCF metric, and the Hopf surface is already counted in the first case, we can assume that the genus $g\geq 1$. In the list we gave in the previous section (Lemma \ref{listtheorem1}), the cases between \ref{elliptic1}-\ref{elliptic4} are elliptic except some 
tori. These are the cases to which the results of this section apply. 
This completes the proof of Theorem \ref{listtheorem2}.


\section{Converse}\label{secconverse}

In this section we will mention the surfaces in the list of Lemma \ref{listtheorem1} that are known to admit LCF metrics. See \cite{structuressdlcf} for a recent survey and \cite{besse} for  references.

\beg{itemize}

\item[Case 1:] 
Among the Hopf surfaces, the primary ones, i.e. the ones homeomorphic to $S^1\times S^3$ admit LCF metrics. The reason is
locally it is a product of a line with a constant curvature space.
Note that if a complex surface is homeomorphic to $S^1\times S^3$ then it is diffeomorphic to it by a result of Kodaira \cite{kodaira}. Among the secondary type Hopf surfaces, the ones obtained by $\mbbz_p$-action on the second component, result in a lens space product $S^1\times L(p,q)$, hence they admit an LCF metric because locally it is as the previous case.

\item[Case 2:] Among the ruled surfaces mentioned, the trivial products are LCF. The product metric on  $\mbb{CP}_1\times \Sigma_g$ admits LCF metric for $g\geq 2$. 

\item[Case 5:] All tori admit flat metrics.

\end{itemize}

\section{The Hermitian case}\label{sechermitianLCF}
In this section we consider the \emph{Hermitian} locally conformally flat structures on complex surfaces. This means that we have an additional compatibility condition, i.e. $J$-invariance relation $g(JX,JY)=g(X,Y)$ for all
vectors $X,Y$.
This case is analyzed by M. Pontecorvo in \cite{pontecorvo}, and also stated without proof in \cite{boyer}. We have one of the following three cases.
\begin{enumerate}
\item \label{ponhopf} A Hopf surface, i.e. finitely covered by a complex surface
$\mbb C^{2*}/\mbbz$ (diffeomorphic to $S^1\times S^3$) with its standard metric of \cite{vaismanlck0}.
\item \label{pongeometricallyruled} A flat $\mbb{CP}_1$ bundle over a Riemann surface $\Sigma_g$ of genus $g\geq 2$. Its metric is locally the product of constant $\pm 1$ curvature metrics.
\item \label{pontorihyperelliptic} A complex torus or a hyperelliptic surface with flat metrics.
\end{enumerate}
In this classification the cases \ref{ponhopf},\ref{pongeometricallyruled} and \ref{pontorihyperelliptic} fall into the cases of $\ref{ourVII}$,
$\ref{ourruled}$ and $\ref{ourelliptic}$ of
our Theorem \ref{listtheorem2} respectively.

\vspace{.5cm}

\vspace{.05in}


{\small
\beg{flushleft}
\textsc{Tuncel\'{\i} \" Un\' ivers\'ites\'i, Mekatron\'{\i}k M\" uhend\' isl\' i\v g\' i, Aktuluk, 
62000, Tuncel\'{i}, T\" urk\'{i}ye.}\\
\textit{E-mail address:} \texttt{\textbf{mkalafat@\,tunceli.edu.tr}}
\end{flushleft}
}


{\small
\beg{flushleft}
\textsc{
Vanderbilt University, Department of Mathematics, Nashville, TN 37240, USA.}\\
\textit{E-mail address:} \texttt{\textbf{caner.koca@\,vanderbilt.edu}}
\end{flushleft}
}

\end{document}